\documentclass[aos,MSNbibl,seceqn,nameyear,dvips]{arximspdf}
\usepackage{graphicx}

\doi{10.1214/11-AOS956}
\volume{40}
\issue{1}
\pubyear{2012}
\firstpage{188}
\lastpage{205}

\makeatletter

\newcommand{\cal}{\mathcal}

\newcommand{\ww}{\omega}

\renewcommand{\overline}{\bar}

\newtheorem{theorem}{Theorem}[section]
\newproclaim{example}{Example}[section]

\makeatother

\begin{document}
\begin{frontmatter}

\title{$T$-optimal designs for discrimination between two polynomial models\thanksref{T1}}
\runtitle{$T$-optimal designs for discrimination}

\thankstext{T1}{Supported in part by the Collaborative Research Center ``Statistical modeling
of nonlinear dynamic processes'' (SFB 823, Teilprojekt C2) of the German Research Foundation
(DFG).}

\begin{aug}
\author[A]{\fnms{Holger} \snm{Dette}\corref{}\ead[label=e1]{holger.dette@rub.de}},
\author[B]{\fnms{Viatcheslav B.} \snm{Melas}\thanksref{t2}\ead[label=e2]{vbmelas@post.ru}}
\and
\author[B]{\fnms{Petr} \snm{Shpilev}\thanksref{t2}\ead[label=e3]{pitshp@hotmail.com}}
\runauthor{H. Dette, V. B. Melas and P. Shpilev}
\affiliation{Ruhr University at Bochum, St. Petersburg State University
and~St.~Petersburg State University}
\address[A]{H. Dette\\
Fakult\"at f\"ur Mathematik \\
Ruhr-Universit\"at Bochum \\
44780 Bochum\\
Germany \\
\printead{e1}}
\address[B]{V. B. Melas \\
P. Shpilev\\
Department of Mathematics \\
St. Petersburg State University \\
St. Petersburg\\
Russia \\
\printead{e2}\\
\hphantom{E-mail: }\printead*{e3}} 
\end{aug}

\thankstext{t2}{Supported in part by Russian Foundation of Basic Research
(Project 09-01-00508).}

\received{\smonth{7} \syear{2011}}
\revised{\smonth{10} \syear{2011}}

%
\vspace*{12pt}
\begin{abstract}
This paper is devoted to the explicit construction of optimal designs
for discrimination between two polynomial regression models of degree
$n-2$ and $n$. In a fundamental paper, Atkinson and Fedorov
[\textit{Biometrika} \textbf{62} (1975a) 57--70] proposed the
$T$-optimality criterion for this purpose. Recently, Atkinson
[\textit{MODA 9}, \textit{Advances in Model-Oriented Design and
Analysis} (2010) \mbox{9--16}] determined $T$-optimal designs for polynomials
up to degree 6 numerically and based on these results he conjectured
that the support points of the optimal design are cosines of the angles
that divide half of the circle into equal parts if the coefficient of
$x^{n-1}$ in the polynomial of larger degree vanishes. In the present
paper we give a strong justification of the conjecture and determine
all $T$-optimal designs explicitly for any degree $n \in\mathbb{N}$. In
particular, we show that there exists a one-dimensional class of
$T$-optimal designs. Moreover, we also present a generalization to the
case when the ratio between the coefficients of $x^{n-1}$ and $x^n$ is
smaller than a certain critical value. Because of the complexity of the
optimization problem, $T$-optimal designs have only been determined
numerically so far, and this paper provides the first explicit solution
of the $T$-optimal design problem since its introduction by
Atkinson and Fedorov
[\textit{Biometrika} \textbf{62} (1975a) 57--70]. Finally, for the
remaining cases (where the ratio of coefficients is larger than the
critical value), we propose a numerical procedure to calculate the
$T$-optimal designs. The results are also illustrated in an example.
\end{abstract}

%

\end{frontmatter}

\section{Introduction}\label{intro}

The problem of identifying an appropriate model in a~class of competing
regression models is of fundamental importance in regression analysis,
and it occurs often in real experimental studies.
It is widely accepted nowadays that good experimental designs can
improve the performance of discrimination, and several authors have
addressed the problem of constructing optimal designs for this purpose;
see \citet{hunrei1965}, \citet{stigler1971},
Atkinson and Fedorov (\citeyear{atkfed1975a,atkfed1975b}),
\citet{hill1978}, \citet{fedorov1981}, \citet
{denfedkha1981}, \citet{studden1982}, \citet{fedkha1986}, \citet
{spruill1990}, Dette (\citeyear{dette1994a,dette1995b}), \citet{dethal1998}, \citet
{songwong1999}, \citet{ucibog2005},
Wiens (\citeyear{wiens2009,wiens2010}) among
many others.
In a fundamental paper, \citet{atkfed1975a}
introduced the $T$-optimality criterion for discriminating between two
competing regression models. As an
example, these authors constructed $T$-optimal designs for a constant
and a quadratic model.
Since its introduction, the problem of determining $T$-optimal designs
has been considered by numerous authors; see \citet{atkfed1975b}, \citet
{ucibog2005}, \citet{wiens2009}, \citet{tomlop2010}, among others.
In order to demonstrate the benefits of the $T$-optimal design, we
display, in Table \ref{tabex}, the simulated power of the $F$-test for
the hypothesis $H_0\dvtx\theta_{2,2} = \theta_{2,3}=0$ in the cubic
regression model
$\eta(x,\theta)=\sum^3_{j=0} \theta_{2j}x^j$ on the interval $[-1,1]$
(with standard normal distributed errors), where observations are taken
according to two designs. The first design is the commonly used
equidistant design with $12$ observations at the four points, $-1,
-1/3, 1/3$ and $1$, respectively, while the second design is a
$T$-optimal design, as considered in this paper, with $8$ observations
at the two points, $-1,1$ and $16$ observations at the two points
$-1/2, 1/2$, respectively. We observe clear advantages (with respect to
the power of the $F$-test) for the $T$-optimal design.

%
\begin{table}[b]
\tablewidth=270pt
\caption{Simulated power of the $F$-test in a cubic
regression model $\sum^3_{j=0}\theta_{2j}x^j$ for the hypothesis of
linear regression model for various values of $\theta_{2,3}$ and
different designs $(\theta_{2,2}=0)$}
\label{tabex}
\begin{tabular*}{\tablewidth}{@{\extracolsep{\fill}}lccccc@{}}
\hline
$\bolds{\theta_{2,3}}$ & $\bolds{0}$ & $\bolds{0.5}$
& $\bolds{1.0}$ & $\bolds{1.5}$ & $\bolds{2.0}$ \\
\hline
$T$-optimal & $0.051$ & $0.104$ & $0.301$ & $0.641$ & $0.896$ \\
Equidistant & $0.053$ & $0.092$ & $0.218$ & $0.438$ & $0.638$ \\
\hline
\end{tabular*}
\end{table}

Since its introduction $T$-optimal designs have found numerous
applications including such important fields as chemistry of
pharmacokinetics; see \citet{atkbogbog1998}, \citet{aspmac2000},
Uci{\'n}s\-ki and Bogacka (\citeyear{ucibog2005}) or \citet{fooduf2011} among others.
The $T$-optimal design problem is essentially a minimax problem, and,
except for very simple models, the corresponding optimal designs are
not easy to find and have to be determined numerically.
In a recent paper,
\citet{dettit2009} discussed the $T$-optimal design problem from a
general point of view and related it to a nonlinear problem in
approximation theory. As an illustration, designs for
discriminating between a linear model and a cubic model without
quadratic term were presented, and it was shown that $T$-optimal
designs are, in general, not unique.

\citet{atkinson2010} considered a similar problem of this
type and studied the problem of discriminating between two competing
polynomial regression models which differ in the degree by two. This
author determined $T$-optimal designs for polynomials up to degree $6$
numerically where the coefficient of~$x^{n-1}$ in the polynomial of
larger degree (say $n$) vanishes. Based on these results
he conjectured
that the support points of the $T$-optimal design are cosines of angles
dividing a half of
circle into equal parts.

The
present paper has two purposes. In particular, we prove the conjecture
raised in \citet{atkinson2010} and
derive explicit solutions of the $T$-optimal design problem for
discriminating between polynomial regression models of degree $n-2$ and
$n$ for any $n \in\mathbb{N}$. Moreover, we also determine the
$T$-optimal designs analytically in
the case when the ratio of the coefficients of the terms~$x^{n-1}$ and~$x^n$ is
sufficiently small.
The situation considered in \citet{atkinson2010} corresponds to the case
where this ratio vanishes, and in this case we show that there exists a
one-dimensional class of $T$-optimal designs.
To the best of our knowledge these results provide the first explicit solution
of the $T$-optimal design problem in a nontrivial situation. Our
results provide further insight into the complicated
structure of the $T$-optimal design problem.
Finally, in the case where the coefficient exceeds the critical value,
we suggest a procedure to determine the $T$-optimal design numerically.

\section{The $T$-optimal design problem revisited}

Consider the classical regression model
%
\begin{equation} \label{10}
y=\eta(x)+\varepsilon,
\end{equation}
where the explanatory variable $x$ varies in the design space $\cal X$,
and observations at different locations, say $x$ and $x'$, are assumed
to be uncorrelated with the same variance. In (\ref{10}) the quantity
$\varepsilon$ denotes a random variable
with mean $0$ and variance $\sigma^2$, and $\eta$ is a function,
which is
called regression function in the literature.
We assume that the experimenter has two parametric models for this
function in mind, that is,
%
\begin{equation}\label{11}
\eta_1(x,\theta_1) \quad\mbox{and}\quad \eta_2(x,\theta_2),
\end{equation}
and the first goal of the experiment is to discriminate between these
two models.
In (\ref{11}) the quantities $\theta_1$ and $\theta_2$ denote unknown
parameters which vary in compact parameter spaces, say $\Theta_1
\subset\mathbb{R}^{m_1}$ and $\Theta_2 \subset\mathbb{R}^{m_2}$, and
have to be estimated from the data.
In order to find ``good'' designs for discriminating between the models
$\eta_1$ and $\eta_2$, we consider approximate designs in the
sense of \citet{Kiefer1974},
which are defined as probability measures on the design space $\mathcal
{X}$ with finite support.
The support points of an (approximate) design $\xi$ give the locations
where observations
are taken, while the weights give the corresponding relative
proportions of total observations to be taken at these points.
If the design $\xi$ has masses $\omega_i>0 $ at the different points
$x_i$ $(i =
1, \ldots, k)$, and $N$ observations can be made by the experimenter,
the quantities
$\omega_i N$ are rounded to integers, say $n_i$, satisfying $\sum
^k_{i=1} n_i =N$, and
the experimenter takes $n_i$ observations at each location $x_i$
$(i=1, \ldots, k)$.

To determine a good design for discriminating between the models~$\eta
_1$ and~$\eta_2$ [\citet{atkfed1975a}] proposed in a fundamental paper to
fix one model, say~$\eta_1$ (more precisely its corresponding parameter
$\theta_1$), and to determine
the design which maximizes the minimal deviation
between the model $\eta_1$ and the class of models defined by $\eta
_2$, that is,
\[
\xi^*=\mathop{\arg\max}_{\xi}\int_{\chi}\bigl(\eta_1(x,\theta_1)-\eta_2
(x,\theta_2^*)\bigr)^{2}\xi(d x),
\]
where the parameter ${\theta}_2^*$ minimizes the expression
\[
{\theta}_2^*=\mathop{\arg\min}_{\theta_2 \in\Theta_2}\int_{\chi}
\bigl(\eta_1(x,\theta_1)-\eta_2(x,\theta_2)\bigr)^{2}\xi(d x).
\]

Note that $\theta^*_2$ is not an estimate, but it corresponds to the
best approximation of the ``given'' model
$\eta_1 (\cdot,\theta_1)$ by models of the form $\{\eta_2(\cdot,
\theta
_2) \mid\theta_2 \in\Theta_2\}$ with
respect to a weighted $L_2$-norm.
Since its introduction the $T$-optimal design problem has found
considerable interest in the literature, and we refer the interested
reader to the work of \citet{ucibog2005} or \citet{dettit2009}, among
others. In general, the determination of $T$-optimal designs is a very
difficult problem, and explicit solutions are---to our best
knowledge---not available except for very simple models with a few parameters.
In this paper we present analytical results for $T$-optimal designs, if
the interest is in the discrimination between two polynomial models
which differ in the degree by two. To be precise,
we consider the case where the regression functions
$\eta_1(x,\theta_1)$ and $\eta_2(x,\theta_2)$ are given by
%
\begin{equation} \label{eta1}\quad
\eta_1(x,\theta_1)=\theta_{1,0}+\theta_{1,1}x+\cdots+\theta
_{1,n-2}x^{n-2}+ \theta_{1,n-1}
x^{n-1}+ \theta_{1,n}x^n
\end{equation}
and
%
\begin{equation} \label{eta2}
\eta_2(x,\theta_2)=\theta_{2,0}+\theta_{2,1} x+\cdots+ \theta
_{2,n-2}x^{n-2},
\end{equation}
respectively, and the design space is given by $\mathcal{X}=[-1,1]$. In
model (\ref{eta1}) the parameter $\theta_1$ is given by
$\theta_1= (\theta_{1,0},\theta_{1,1},\ldots, \theta_{1,n-2}, b
\theta
_{1,n}, \theta_{1,n})^T$, where the ratio of the coefficients
corresponding to the highest powers $b=\theta_{1,n-1}/\theta_{1,n}$ and
the parameter $\theta_{1,n}$ specify the deviation from a polynomial of
degree $n-2$.\vadjust{\goodbreak}

In the following discussion, we define
%
\begin{eqnarray} \label{etaq}
\overline{\eta}(x,\alpha,b,\theta_{1,n}) &=& \eta_1(x,\theta
_1)-\eta
_2(x,\theta_2)
\nonumber\\[-8pt]\\[-8pt]
&=& \alpha_0+\alpha_1
x+\cdots+\alpha_{n-2}x^{n-2} +\theta_{1,n} (bx^{n-1}+x^n),\nonumber
\end{eqnarray}
where we use the notation
$\alpha_i= \theta_{1,i}-\theta_{2,i}$ $(i=0,\ldots,n-2)$;
then the problem of finding the $T$-optimal design for the models $\eta
_1$ and $\eta_2$ can be reduced to
\[
\xi^*=\mathop{\arg\max}_{\xi}\int_{\chi}\bigl(\alpha_0^* +\alpha
_1^*x+\cdots
+\alpha_{n-2}^*x^{n-2}+\theta_{1,n}
(bx^{n-1}+x^n)\bigr)^{2}\xi(d x),
\]
where $\alpha^*=(\alpha_1^*,\ldots,\alpha_{n-2}^*)^T$ is a
vector minimizing the expression
\[
\alpha^*=\mathop{\arg
\min}_{\alpha}\int_{\chi}(\overline{\eta}(x,\alpha,b,\theta
_{1,n}))^{2}\xi(d
x).
\]

It is now easy to see that for a fixed value of $b=\theta
_{1,n-1}/\theta
_{1,n}$, the $T$-optimal design does not depend on the parameter
$\theta_{1n}$.
In the next section we give the complete solution of the $T$-optimal
design problem if the absolute value of the parameter $b=\theta
_{1,n-1}/\theta_{1,n}$ less or equal to
some critical value.

\section{\texorpdfstring{$T$-optimal designs for small values of $|b|=|\theta_{1,n-1}/\theta_{1,n}|$}
{T-optimal designs for small values of |b|=|theta_{1,n-1}/theta_{1,n}|}}
\label{sec3}

Throughout this section we assume that the parameter $b$ satisfies
%
\begin{equation} \label{ass1}
|b|\,{=}\,|\theta_{1,n-1}/\theta_{1,n}|\,{\leq}\,n \biggl(1\,{-}\,\cos
\biggl(\frac{\pi}{n}\biggr)\biggr)\bigg/\biggl(1\,{+}\,\cos
\biggl(\frac{\pi}{n}\biggr)\biggr)\,{=}\,n\tan^2\biggl(\frac{\pi}{2n}\biggr);\hspace*{-40pt}
\end{equation}
then it is easy to see that all points
%
\begin{equation} \label{supp}
t_i^*(b)=-\biggl(1+\frac{|b|}{n}\biggr)\cos\biggl(\frac{i\pi
}{n}\biggr)-\frac{|b|}{n},\qquad
i=1,\ldots,n,
\end{equation}
are located in the interval $[-1,1]$.
Our first result gives an explicit solution of the $T$-optimal design
problem in the case $b=\theta_{1,n-1}=0$ and---as a by-product---proves
the conjecture raised in \citet{atkinson2010}.
%
\begin{theorem}\label{theo31}
A design $\xi$ is $T$-optimal for discriminating between models
(\ref{eta1}) and (\ref{eta2}) with $\theta_{1n-1}=0$ on the interval $[-1,1]$
if and only if it
can be represented in the form $\xi=(1-\alpha)\xi_1 + \alpha\xi_2$,
where $\alpha\in[0,1]$, the measures~$\xi_1$ and $\xi_2$ are defined by
%
\begin{equation} \label{xi12}\qquad
\xi_1=\pmatrix{
t^*_{1}(0) & \cdots& t^*_{n}(0)\cr
\ww^*_1 & \cdots& \ww^*_{n}},\qquad
\xi_2=\pmatrix{
- t^*_{n}(0) & \cdots& -t^*_{1}(0)\cr
\ww^*_n & \cdots& \ww^*_{1}}
\end{equation}
and the weights and support points are given
by
%
\begin{eqnarray} \label{weight}\qquad
\ww^*_i&=&\frac{2}{n}\sin^2\biggl(\frac{i\pi}{2n}\biggr),\qquad
\ww^*_{n-i}=\frac{2}{n}\cos^2\biggl(\frac{i\pi}{2n}\biggr),\qquad
i=1,\ldots,\biggl\lfloor\frac{n}{2}\biggr\rfloor,\nonumber\\[-8pt]\\[-8pt]
\ww^*_n&=&\frac{1}{n}\nonumber
\end{eqnarray}
and (\ref{supp}) for $b=0$, respectively.\vadjust{\goodbreak}
\end{theorem}
\begin{pf}
It was proved by \citet
{dettit2009} (see Theorem 2.1) that any $T$-optimal design on the
interval $[-1,1]$ for discriminating between the polynomials $\sum
^{n-2}_{j=0} \theta_{2,j}x^j$ and
\[
\eta_1(x,\theta_1)= \sum^{n-2}_{j=0} \theta_{1,j} x^{j} + \theta_{1n}x^n
\]
(note that $\theta_{1n-1}=0$) is supported at the set of the extremal points
\[
\mathcal{A} = \Bigl\{ x \in[-1,1] \bigm| \psi^*(x) = \sup_{t
\in
[-1,1]} |\psi^*(t)| \Bigr\},
\]
where $\psi^*(x)=\eta_1(x,\theta_1) - \sum^{n-2}_{j=0} \overline
\theta
_{2j} x^j$ and
%
\begin{equation} \label{appr1}\qquad
\overline\theta_2 = (\overline\theta_{2,0},\ldots,\overline\theta
_{2,n-2})^T = \mathop{\arg\min}_{\theta_2 \in\mathbb{R}^{n-1}} \sup_{x \in
[-1,1]} \Biggl|\eta_1(x,\theta_1) - \sum^{n-2}_{j=0} \theta_{2,j}x^j
\Biggr|
\end{equation}
is the parameter corresponding to the best approximation of $\eta
_1(x,\theta_1)$ with respect to the sup-norm. By a standard result in
approximation theory [see \citet{achiezer1956}, Sections 35 and 43] it
follows that the solution of the problem (\ref{appr1}) is unique and
given by $\psi^*(x)=\theta_{1,n}2^{-(n-1)} T_n(x)$, where
$T_n(x)=\cos(n\arccos x)$ is the $n$th Chebyshev polynomial of the first
kind. Note that $T_n(x)$ is an even or odd polynomial of degree $n$
with leading coefficient $2^{n-1}$ [see \citet{szego1975}]. The
corresponding extremal points are given by
$x_0= t^*_1(0)=-1$, $x_i= t^*_i(0)=- \cos\frac{i\pi}{n}$,
$i=1,\ldots
,n-1$, $x_n= t^*_n(0)=1$.

Now it follows from Theorem 2.2 in \citet{dettit2009} that
a design~$\xi^*$ is $T$-optimal if and only if it satisfies the system
of linear equations
%
\begin{equation}\label{appr2}
\int_{\mathcal{A}} \psi^*(x) x^k \,d \xi^*(x)=0,\qquad k=0,\ldots,n-2.
\end{equation}
[Note that in the case of linear models the necessary condition in
Theorem~2.2 in \citet{dettit2009} is also sufficient.]
Therefore for proving that
$\xi^*_1=\xi_1$ is a $T$-optimal design, it is sufficient to verify the
identities
%
\begin{equation} \label{weight1}
\int\psi^* (x) \,d \xi^*_1 (x) = \theta_{1,n} 2^{-(n-1)} (-1)^n
\sum^n_{i=1}(-1)^i x^k_i\omega^*_i =0
\end{equation}
($k=0,1,\ldots,n-2$), which will be done in the \hyperref[app]{Appendix}.
In a similar way we can check that the design $ \xi_2^*$ in (\ref
{xi12}) is a $T$-optimal design. Note that
\[
\operatorname{supp} (\xi^*_1)
\cup\operatorname{supp} (\xi^*_2) =
\biggl\{ x_i= - \cos\biggl( \frac{\pi}{n} i \biggr) \Bigm|  i=0,\ldots,n
\biggr\} = \mathcal{A},
\]
because $t^*_{n-i}(0)=-t^*_i(0)$.
Moreover, (\ref{appr2}) defines a system of linear equations of the
form $F \omega=0$ for the vector $\omega=(\omega_0,\ldots,\omega_n)^T$
of the $T$-optimal design~$\xi^*$, where the matrix $F$ is given by
$F=((-1)^i x^k_i)^{k=0,\ldots,n-2}_{i=0,\ldots,n} \in\mathbb{R}^{n-1
\times n+1}$
and has rank $n-1$.
Additionally, the components of the vector $\omega$ satisfy $\sum
^n_{i=0} \omega_i=1$. Therefore the set of solutions has dimension $1$.
Because the vectors of weights corresponding to the designs $\xi^*_1$
and $\xi^*_2$ are given by $\omega^{(1)}=(0,\omega^*_1,\ldots,\omega
^*_n)^T$ and
$\omega^{(2)}=(\omega^*_n,\ldots,\omega^*_1,0)^T$ and are therefore
linearly independent (note that $\omega^*_i >0, i=1,\ldots,n$), any
vector of weights corresponding to a $T$-optimal design must be a
convex combination of~$\omega^{(1)}$ and~$\omega^{(2)}$. Consequently,
any $T$-optimal design can be represented in the form $\xi=(1-\alpha)
\xi^*_1 + \alpha\xi^*_2$, which proves the assertion of Theorem \ref
{theo31}.
\end{pf}

Note that the $T$-optimal design is not unique in the case $b=0$. On
the other hand, the $T$-optimal
designs are unique, whenever $\theta_{1,n-1} \not=0$, and, if the
ratio $| \theta_{1,n-1} / \theta_{1,n}|$
is not too large, the $T$-optimal designs can also be found explicitly
as demonstrated in our following
result.
%
\begin{theorem}\label{theo32}
If the parameter $b=\theta_{1,n-1}/\theta_{1,n}$ satisfies (\ref
{ass1}), then there exists a unique $T$-optimal design on the interval
$[-1,1]$ for discriminating between models (\ref{eta1}) and (\ref
{eta2}). For positive $b$ this design has
the form
%
\begin{equation} \label{des1}
\xi^*=\pmatrix{
t^*_{1}(b) & \cdots& t^*_{n}(b)\cr
\ww^*_1 & \cdots& \ww^*_{n}},
\end{equation}
where the points $t^*_i(b)$ and weights $w^*_i(b)$ are defined in
(\ref{supp}) and (\ref{weight}), respectively [note that
$t^*_1(b)\geq
-1,t^*_n(b)=1$].
The $T$-optimal design for negative~$b$ has the form
\[
\xi^*=\pmatrix{
-t^*_{n}(b) & \cdots& -t^*_{1}(b)\cr
\ww^*_n & \cdots& \ww^*_{1}}
\]
[note that $-t^*_n(b)= -1,-t^*_1(b)\leq1$].
\end{theorem}
\begin{pf}
We consider the case $0< b
\leq n (1 - \cos(\frac{\pi}{n}))/(1 +
\cos(\frac{\pi}{n}))$ where
direct calculations show that the points $t^*_i(b),
i=1,\ldots,n$, are contained in the interval $[-1,1]$. Moreover, these
points are the
extremal points of the polynomial
%
\begin{equation} \label{exp}
c_nT_n\biggl(\frac{-x-{b/n}}{ 1+ {b/n} }\biggr),\qquad c_n = (-1)^n
\biggl(\frac{1}{2}\biggr)^{n-1}\biggl(1+ \frac{b}{n}\biggr)^{n},
\end{equation}
where $T_n$ is the Chebyshev polynomial of the first kind.
For later purposes we note that the coefficient of $x^{n-1}$ in this
polynomial is equal to
%
\begin{equation} \label{erg}
\sum_{i=1}^n
\biggl[\biggl(1+\frac{b}{n}\biggr)u_i+\frac{b}{n}\biggr]= b,
\end{equation}
where $u_1,\ldots,u_n$ are the roots of the polynomial $T_n(x)$, that
is, $u_i=\cos(\frac{2i-1}{2n}\pi) $
$(i=1,\ldots,n)$, $\sum_{i=1}^nu_i=0$.
It can be shown by a standard argument in approximation theory [see
\citet{achiezer1956}, Sections 35 and 43] that $\theta_{1n}\psi^*(x)$ with
\[
\psi^*(x)= c_n T_n \biggl( \frac{-x-{b/n}}{1+{b/n}}
\biggr)
\]
is the unique solution of the extremal problem
\[
\min_{\theta_2 \in\mathbb{R}^{n-1}} \sup_{x \in[-1,1]} \Biggl|\eta
_1(x,\theta_1) - \sum^{n-2}_{j=0} \theta_{2,j}x^j \Biggr|,
\]
where $ \eta_1 (x,\theta_1)=\sum^n_{j=0}\theta_{1,j}x^j$. Therefore by
Theorems 2.1 and 2.2 in \citet{dettit2009}, a $T$-optimal design is
supported at
the $n$ extremal points $t_1^*(b),\ldots,t_n^*(b)$
[note that we use $b \leq n \tan^2 (\frac{\pi}{2n})$ at this point,
which implies $|t^*_j(b)|\leq1; j=1,\ldots,n$]
and the weights are determined by (\ref{appr2}). Because the set of
extremal points is given by
$\mathcal{A}=\{t^*_1(b),\ldots,t^*_n(b)\}$, this system reduces to
%
\begin{equation} \label{weight2}
\sum^n_{i=1}t^{*k}_i(b)(-1)^i \omega^*_i =0,\qquad k=0,1,\ldots,n-2,
\end{equation}
and we will prove in the \hyperref[app]{Appendix} that the weights given in (\ref
{weight}) define a~solution of (\ref{weight2}). Therefore
the design $\xi^*$ specified in (\ref{des1}) is a
$T$-optimal design for $0 < b\leq n (1 - \cos\pi/n)/(1 +
\cos\pi/n)$. Since the function $\psi^*(x)$ is unique, any $T$-optimal
design is supported at the points
$t_1^*(b),\ldots,t_n^*(b)$ [see Theorem 2.1 in \citet{dettit2009}]. By
Theorem 2.2 in the same reference, it follows that the weights of any
$T$-optimal design satisfy the
system of linear equations (\ref{weight2}) with $\omega^*_i = \omega
_i$ and
$\sum^n_{i=1}\ww_i=1$.
Since $\psi^*(t^*_i(b))=(-1)^i$ $(i=1,\ldots,n)$ we can rewrite this
system as
%
\begin{equation} \label{eq1}
F\ww=e_n,
\end{equation}
where $\ww=(\ww_1,\ldots,\ww_n)^T$ is the vector of weights, the last
row of the matrix~$F$ is given by $(1,\ldots,1)$ and corresponds to the
condition $\sum^n_{i=1}\omega_i=1$, $e_n=(0,\ldots,0,1)^T \in\mathbb
{R}^n$ denotes the $n$th unit vector
and the columns of the matrix~$F$ are given by
\[
a_i=(-1)^i(1,t^*_i(b),\ldots,(t^*_i(b))^{n-2},\psi^*(t^*_i(b)))^T,\qquad
i=1,2,\ldots,n.
\]
The remaining assertion of Theorem \ref{theo32} follows if we prove
that $\det F\ne0$, which implies that the solution of (\ref{eq1}), and
therefore the $T$-optimal design, is unique. For this purpose assume
that the opposite holds. In this case the rows of the matrix $F$ would
be linearly dependent, and there exists a vector
$h=(h_1,\ldots,h_{n-1},1)^T$ such that $a_i^Th=0, i=1,2,\ldots,n$. But
the function $k(x)=(1,x,\ldots ,x^{n-2},\psi^*(x))^Th$ is a polynomial
of degree $n$ with coefficient of $x^{n-1}$ given\vadjust{\goodbreak} by $b$. Since $a_ih
=k(t^*_i(b))=0$ this polynomial has roots at the points~$t_i^*(b)$,
moreover
\[
\sum_{i=1}^nt_i^*(b)= -b
-\sum_{i=1}^n\biggl(1+\frac{b}{n}\biggr)\cos\biggl(\frac{i\pi
}{n}\biggr)=-b
+1+\frac{b}{n}.
\]
However, by (\ref{erg}) the sum of the roots must equal $-b$ by
Vieta's formula.
This contradiction proves that $\det F\ne0$. Therefore the system of
equations in (\ref{eq1}) has a unique solution, which means that the
$T$-optimal design is unique.

The case of negative $b$ is considered in a similar way, and the
details are omitted for the sake of brevity.
\end{pf}

The critical values $b_n^*=n\tan^2(\frac{\pi}{2n})$ for
various values of $n \in\mathbb{N}$ are displayed in Table \ref{tab6}.
Theorems \ref{theo31} and \ref{theo32} give an explicit solution of the
%
\begin{table}
\caption{The critical values
$b^*_n = n \tan^2(\frac{\pi}{2n})$ for
various values $n \in\mathbb{N}$}\label{tab6}
\begin{tabular*}{\tablewidth}{@{\extracolsep{\fill}}lcccccccc@{}}
\hline
$\bolds{n}$ & $\bolds{3}$ & $\bolds{4}$ & $\bolds{5}$ & $\bolds{6}$
& $\bolds{7}$ & $\bolds{8}$ & $\bolds{9}$ & $\bolds{10}$ \\
\hline
$b^*_n$ &1& 0.6864&
0.5280& 0.4306& 0.3646& 0.3168& 0.2801& 0.2509 \\
\hline
\end{tabular*}
\vspace*{-3pt}
\end{table}
$T$-optimal design problem for discriminating between a polynomial
regression of degree $n-2$ and~$n$, whenever $|b| =|\theta
_{1,n-1}|/|\theta_{1,n}| \leq b_n$. In the opposite case the solution
is not so transparent and will be discussed in the following
section.\vspace*{-3pt}


\section{$T$-optimal designs for large values of $|b|$}\label{sec4}

In this section\vspace*{2pt} we consider the case $|b|\geq
n\tan^2(\frac{\pi}{2n})$, for which the $T$-optimal design
cannot be found explicitly.
Therefore we present a numerical method
to determine the optimal designs. The method was described by \citet
{detmelpep2004a}
in the context of determining optimal designs for estimating individual
coefficients in a polynomial regression
model [see also \citet{melas2006}], and for the sake of brevity, we only
explain the basic principle.
For this purpose we rewrite the function
$\overline{\eta}$ in (\ref{etaq}) as
%
\begin{equation} \label{etaq2}
\overline{\eta}(x,\alpha,\bar b)= \alpha_0+\alpha_1
x+\cdots+\alpha_{n-2}x^{n-2} + \theta_{1n-1} (x^{n-1}+\bar b x^n),
\end{equation}
where $\bar b=1/b = \theta_{1n}/\theta_{1n-1}$.
Note that for fixed $\bar b$, the $T$-optimal design is independent of
the parameter
$\theta_{1n-1}$ and that the choice
\[
\bar b \in\biggl[-\frac
{1}{n}\cot^2\biggl(\frac{\pi}{2n}\biggr),\frac
{1}{n}\cot^2\biggl(\frac{\pi}{2n}\biggr)\biggr]
\]
corresponds to the case $|b|\geq n\tan^2(\frac{\pi}{2n})$
considered in this section.
In order to express the dependence on the parameter $\bar b$, we use
the notation $t^*_i(\bar b)$ for the support points and $\omega
^*_i(\bar b)$ for
the weights of the $T$-optimal design in this section.

The main idea of the algorithm
is a representation of the support\break points~$t_i^*(\bar b) $ and
corresponding\vadjust{\goodbreak} weights $\omega_i^*(\bar b) $ in terms of a Taylor
series, where the coefficients can be determined
explicitly as soon as the design is known for a particular point~$\bar
b $. The algorithm proceeds in several steps:
\begin{longlist}[(1)]
\item[(1)] \textit{Initialization}:
In the present situation the point $\bar b$ is given by
$\bar b=0$, which corresponds to the situation of discriminating
between a polynomial of degree \mbox{$n-2$} and $n-1$.
For this case it follows from \citet{dettit2009} that the $T$-optimal
design coincides with
the $D_1$-optimal design. This design\vspace*{1pt} has been determined explicitly
by \citet{studden1980}
and puts masses $\omega_i (0) =\frac{1}{n-1}$ at the points $t_{i} (0)
= \cos(\frac{(i-1)\pi}{n-1})$ ($ i=2,\ldots, n-1$)
and masses $\omega_1 (0) = \omega_{n} (0) = \frac{1}{2(n-1)}$ at the
points $t_{1} (0) =-1$ and $t_{n} (0) =1$.
\item[(2)] \textit{The dual problem}:
For
the constructions of the Taylor expansion we now associate to each vector
\begin{eqnarray*}
\tau\in{\cal U} &=&\Biggl\{(t_2,\ldots,t_{n-1},\omega_1,\ldots
,\omega
_{n-1})^T
\Bigm| -1<t_2<\cdots<t_{n-1}<1;\\
&&\hspace*{176pt}\omega_i>0, \sum^{n-1}_{j=1}\omega_j<1 \Biggr\},
\end{eqnarray*}
a design with $n$ support points defined by
\[
\xi_\tau=\pmatrix{
-1&t_2&\cdots&t_{n-1}&1\cr
\omega_1&\omega_2&\cdots&\omega_{n-1}&\omega_n}.
\]
As pointed out in the previous discussion,
there exists a corresponding extremal problem defined by
%
\begin{equation} \label{extr}
\inf_{q\in\mathbb{R}^{n-1}} \sup_{x\in[-1,1]}|\bar
bx^n+x^{n-1}-\bar
f^T(x)q|
\end{equation}
with a unique solution corresponding to
the $T$-optimal design problem under consideration,
where we use the notation
$\bar f^T(x) = (1,x,\ldots,x^{n-2})$.
\item[(3)] \textit{The necessary condition}:
For each vector $q$ in (\ref{extr}), define vectors $d_q = (q^T,1,\bar
b)^T, \Theta=(q,\tau)$ and
a quadratic form
\[
H(\Theta, \bar b)= H(q,\tau,\bar b)=d_q^TM(\xi_\tau)d_q,
\]
where $M(\xi_\tau)$ is the information matrix of the design
$\xi_\tau$ for the regression model~(\ref{etaq2}).
It then follows by similar results as in \citet{detmelpep2004a} that
the design $\xi_{\tau^*}$ is a $T$-optimal design for
discriminating between the polynomials of degree $n$ and $n-2$,
and the vector $q^*$ is a solution of an extremal problem (\ref{extr})
if the points $\Theta^*=(q^*,\tau^*)\in
\mathbb{R}^{n-1} \times{\cal U} $ are the unique solution of the system
\[
\frac{\partial}{\partial\Theta}H(\Theta,\bar b) \bigg|_{\Theta
=\Theta
^*} =0,
\]
such that the inequality
$
|d_{q^*}^Tf(x)|^2\leq d_{q^*}^T M(\xi_{\tau^*})d_{q^*}
$
holds for all $x\in[-1,1]$.\vspace*{1pt}\vadjust{\goodbreak}
\item[(4)] \textit{Taylor expansion of the optimal solution}:
The function
\[
\Theta^* \dvtx \cases{
I \longrightarrow\mathbb{R}^{3n-4}, \cr
\bar b \longrightarrow\Theta^*(\bar b )=(\Theta_1^*(\bar b ),\ldots
,\Theta_{3n-4}^*(\bar b ))=(q^*(\bar b )^T,\tau^*(\bar b )^T),}
\]
which\vspace*{1pt} maps the parameter $\bar b \in I =
[-\frac{1}{n}\cot^2 (\frac {\pi}{2n}),\frac
{1}{n}\cot^2(\frac{\pi}{2n})]$ to the coordinates of the best
approximation $q^*(\bar b )$ and the support points $t^*_i(\bar b)$ and
weights $\omega^*(\bar b )$ of the $T$-optimal design, is a real
analytical function. The coefficients in the corresponding Taylor
expansion,
\[
\Theta^*(\bar b)=\Theta^*(\bar b_0)+\sum_{j=1}^\infty\Theta ^*(j,\bar
b_0)(\bar b-\bar b_0)^j
\]
in a neighborhood of any point $\bar b_0\in I$, can be calculated by
the recursive formulas
\begin{eqnarray}
\Theta^*(s+1,\bar b_0)= -\frac{1}{(s+1)!} J^{-1}(\bar b_0)
\biggl(\frac
{d}{db}\biggr)^{s+1} g\bigl(\Theta^*_{(s)}(\bar b),\bar b\bigr)\bigg|_{\bar b =
\bar b_0},\nonumber\\
&&\eqntext{s=0,1,2,\ldots,}
\end{eqnarray}
where
\begin{eqnarray*}
\Theta^*_{(s)}(\bar b)&=&\Theta^*_{(s)}(\bar
b_0)+\sum_{j=1}^s\Theta^*(j,\bar b_0)(\bar b-\bar b_0)^j,
\\
g(\Theta,\bar b) &=&\frac{\partial}{\partial\Theta} H(\Theta,\bar
b), \\
J(\bar
b_0)&=&\biggl(\frac{\partial^2}{\partial\Theta_i\,\partial\Theta_j}
H(\Theta,\bar b)\biggr) \bigg|_{\Theta=\Theta^*(\bar b_0)}.
\end{eqnarray*}
\end{longlist}
We can use this procedure to calculate the
$T$-optimal design for discriminating between polynomials of degree $n$
and $n-2$
in the cases which are not covered by Theorems \ref{theo31} and \ref
{theo32}. We illustrate the methodology
in the following example.
%
\begin{example}
Consider the $T$-optimal design problem for a model of degree~$5$ and a
cubic polynomial model. Note that for $n=5$, we have
$n\tan^2(\frac{\pi}{2n}) \simeq 0.528$. Therefore if $b\in[0,0.528]$, a
$T$-optimal design is given by Theorem \ref{theo31}, that is,
\begin{eqnarray*}
\xi^*_{T}&=&\pmatrix{
t_1(b) & t_2(b) & t_3(b) & t_4(b) & 1\cr
0.038 & 0.138 & 0.262 & 0.362 & \frac{1}{5}},
\\
t^*_i(b)&=&-\biggl(1+\frac{b}{5}\biggr)\cos\biggl(\frac{i\pi
}{5}
\biggr)-\frac{b}{5},\qquad
i=1,\ldots,5.
\end{eqnarray*}
In order to construct the $T$-optimal design on the interval
$[0.528,\infty]$,
we introduce the notation $\overline b=1/b \in[0, 1.894]$.
With\vspace*{1pt} the results of the previous paragraph we obtain a Taylor expansion
for the interior support points
$t_2^*(\overline{b}),t_3^*(\overline{b}),t_4^*(\overline{b})$ and
weights $\ww_1^*(\overline{b}),
\ww_2^*(\overline{b}),\ww_3^*(\overline{b}),\ww_4^*(\overline{b})$
of the $T$-optimal design for discriminating between a cubic and a
polynomial of degree $5$ where
$\overline b= \theta_{1n}/\theta_{1n-1}$. By the results of \citet
{studden1980}, the vector of support points and weights corresponding
to the center of the expansion at the point $\bar b_0=0$
is explicitly known; that is,
\[
(t_2^*(0), t_3^*(0), t_4^*(0), \omega_1^*(0),\ldots,\omega_4^*(0)) =
\bigl(-\tfrac{1}{\sqrt{2}},0,
\tfrac{1}{\sqrt{2}},\tfrac{1}{8},\tfrac{1}{4},\tfrac{1}{4},
\tfrac{1}{4}\bigr).
\]
At the first step we use a Taylor expansion at the point $\bar b_0=0$
to determine the $T$-optimal design
for $\bar b \in[0,0.4]$. When we have found the vector $\Theta^*(0.4)$
we construct a further Taylor expansion at the point $\bar b_0=0.4$,
%
\begin{figure}

\includegraphics{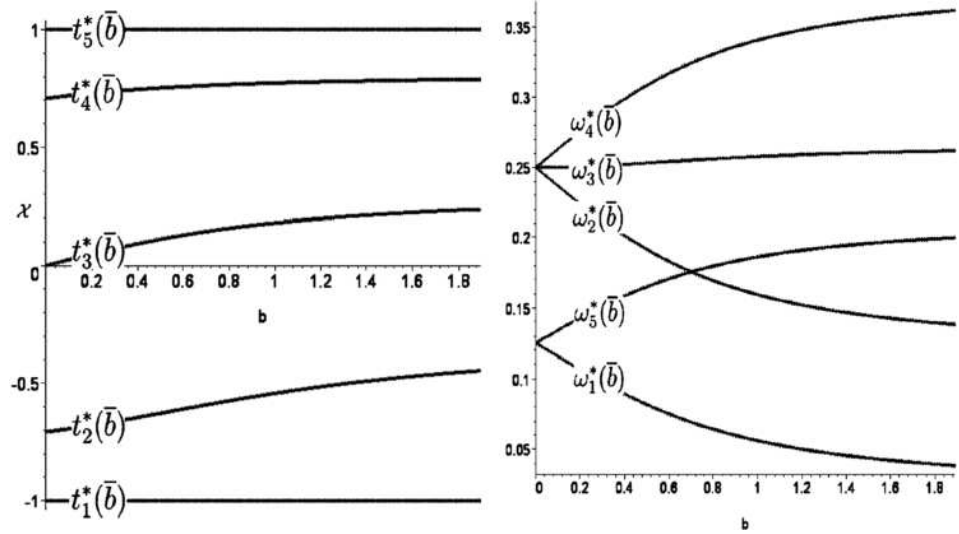}

\caption{The support points (left panel) and weights (right panel) of
the $T$-optimal design for discriminating between a polynomial of
degree $3$ and $5$ for various values of $\overline{b} =1/b \in
[0,1.894]$.}
\label{fig3}
\end{figure}
and this process is continued in order to determine the vector $\Theta
^*(\overline{b})$
for any value $\overline{b}\in[0,1.894]$. The support points and
weights are depicted in Figure \ref{fig3}
as a~function of the parameter $\bar b=1/b=\theta_{1n}/\theta_{1n-1}$.
Note that in all cases $b \neq0$ the $T$-optimal design for
discriminating between a polynomial of degree $5$ and $3$ is supported
at $5$ points.
\end{example}

\section{Concluding remarks and further discussion}

In this paper we have determined $T$-optimal designs for discriminating
between two rival polynomial regression models of degree $n-2$ and $n$.
To the best of our knowledge these results provide the first analytic solution
of a $T$-optimal discriminating design problem with an arbitrary number of
parameters in the regression model.

It should be pointed out that the results depend on the ratio of the
coefficients of the terms $x^{n-1}$ and $x^n$ in the polynomial of
larger degree, which is a well-known feature of the $T$-optimality
criterion. Therefore the designs derived here are local in the sense of
\citet{chernoff1953}. Usually locally optimal designs serve as a
benchmark for commonly used designs, as demonstrated in the example in
the \hyperref[intro]{Introduction}. Moreover, locally optimal designs form the basis for
more sophisticated design strategies, which require less knowledge
about the model parameters such as Bayesian or standardized maximin
optimality criteria [see \citet{chalverd1995} or \citet{dette1997}, among
others]. This extension was already mentioned in the pioneering work of
Atkinson and Fedorov (\citeyear{atkfed1975a,atkfed1975b}) and we
conclude this paper with a~brief discussion of a first explicit result
on maximin $T$-optimal designs for the polynomial regression models.

To be precise, consider the situation where the ratio $b= \theta
_{1,n-1}/\theta_{1,n}$ cannot be exactly specified, but prior knowledge
suggests that $b \in I$ for some interval $I \subset\mathbb{R}$.
Without loss of generality, assume $\theta_{1,n}=1$; then following
\citet{atkfed1975a}, a maximin optimal discriminating design maximizes
the expression
%
\begin{equation}\label{maxi}
\inf_{b \in I} \inf_{\theta_2 \in\mathbb{R}^{n-1}} \int^1_{-1}
\Biggl(x^n +
bx^{n-1} + \sum^{n-2}_{j=0} \theta_{2,j}x^j\Biggr)^2 \,d\xi(x).
\end{equation}
The following result provides a solution of this optimal design problem
for specific intervals $I \subset\mathbb{R}$.
%
\begin{theorem}\label{thm51}
\textup{(a)} If $I=\mathbb{R}$, the maximin $T$-optimal discriminating design is
given by
%
\begin{equation}\label{maxi1}
\xi^*_{MM} = \pmatrix{
t^*_0 & t^*_1 & \cdots& t^*_{n-1} & t^*_n \vspace*{2pt}\cr
\dfrac{1}{2n} & \dfrac{1}{n} & \cdots& \dfrac{1}{n} & \dfrac{1}{2n}},
\end{equation}
where the support points are defined by
\[
t^*_i = \cos\biggl(\frac{n-i}{n} \pi\biggr),\qquad i=0,\ldots,n.
\]

\textup{(b)} Assume that $I=(-\infty, -b_0]$ or $I=[b_0, \infty)$. If $b_0 \ge
0$, then the maximin $T$-optimal discriminating design coincides with
the $T$-optimal discriminating design determined in Sections \ref{sec3} and
\ref{sec4}
for the value $b= b_0$.

In particular, if $b_0=0$, then all designs specified in Theorem \ref
{theo31} are maximin $T$-optimal discriminating designs.
\end{theorem}
\begin{pf}
In order to prove part (a), note
that for $I=\mathbb{R}$ criterion (\ref{maxi}) reduces to
\[
\sup_\xi\inf_{\theta\in\mathbb{R}^{n}} \int^1_{-1} \Biggl(x^n + \sum
^{n-1}_{j=0} \theta_{2,n-1} x^j\Biggr)^2 \,d \xi(x),
\]
which corresponds to the $T$-optimal discriminating design problem for
a~polynomial of degree $n$ and $n-1$. By the results in \citet
{dettit2009}, the solution of this problem coincides with the
$D_1$-optimal design, which is given by (\ref{maxi1}) [see \citet
{studden1980}].

For a proof of part (b), observe that
\begin{eqnarray*}
&& \sup_\xi\inf_{b \in I} \sup_{\theta_2 \in\mathbb{R}^{n-1}}
\int
^1_{-1} \Biggl(x^n + bx^{n-1} + \sum^{n-2}_{j=0} \theta_{2,j}
x^j
\Biggr)^2 \,d \xi(x) \\
&&\qquad = \inf_{b \in I} \sup_\xi\sup_{\theta_2 \in
\mathbb
{R}^{n-1}} \int^1_{-1} \Biggl(x^n + bx^{n-1} + \sum^{n-2}_{j=0}
\theta
_{2,j} x^j\Biggr)^2 \,d \xi(x) =: \inf_{b \in I} R(b),
\end{eqnarray*}
where the last equality defines the function $R$ in an obvious manner.
We now consider the case $I=[b_0, \infty)$ with $b_0 \geq0$ and show
that the function $R$ is increasing on $\mathbb{R}^+$, which implies
%
\begin{equation}\label{maxi2}
\inf_{b \in I} R(b) = R(b_0)
\end{equation}
and proves the assertion for the case $I=[b_0, \infty)$. Recall the
definition of $b^*=n$ $\tan^2(\pi/2n)$ in (\ref{ass1}); then the proof
of Theorem \ref{theo31} shows that
for $\in(0,b^*]$
\[
R(b)= \biggl(1+ \frac{b}{n}\biggr)^{2n} \frac{1}{2^{2n-2}},
\]
which is obviously increasing with respect to the argument $b$. If $R$
would be not increasing on the remaining region $\mathbb{R}^+
\setminus(0,b^*]$, then there would exist real numbers $b_2 > b_1 >
b^*$, such that $R(b_1)= R(b_2)$ with corresponding extremal polynomials
\[
L_i(x)=x^n + b_ix^{n-1} + q^T_i \bar f(x),\qquad i=1,2,
\]
where $\bar f(x)=(1,x,\ldots,x^{n-2})^T$ and
\[
q_i = \mathop{\arg\min}_{q \in\mathbb{R}^{n-1}} \int^1_{-1} \bigl(x^n + b_i
x^{n-1} + q^T \bar f(x)\bigr)^2 \,d \xi(x).
\]
This yields
\[
{\sup_{x \in[-1,1]}} | L_1(x) |= {\sup_{x \in[-1,1]} }|L_2(x)
|= \sqrt{R(b_1)} = \sqrt{R(b_2)}.
\]
By the discussion in Section \ref{sec4}, the polynomials $L_1, L_2$ can be
chosen such that they coincide at the boundary points\vadjust{\goodbreak} of the interval
$[-1,1]$ (note that for $b > b^*$ the support of
the optimal discriminating design always contains
both boundary points $-1$ and $1$). Therefore a simple argument shows
that there exist $n-2$ other points in the interior of the interval
$(-1,1)$, where the polynomials must coincide. Consequently,
$L_1(\tilde t_j) = L_2 (\tilde t_j)$ for $n$ points $\tilde t_1, \ldots,
\tilde t_n \in[-1,1]$, which shows that the polynomials are identical.
This yields $b_1= b_2$, and because of this contradiction the
monotonicity of the function $R$ has been established, which proves
(\ref{maxi2}) and part (b) in the case $I=[b_0,\infty)$. The remaining
case $I=(-\infty,-b_0]$ can be proved by similar arguments, and the
details are omitted for the sake of brevity.
\end{pf}

Theorem \ref{thm51} provides the solution to maximin $T$-optimal
discriminating design problems
for specific intervals $I\subset\mathbb{R}$. In particular, it
identifies the worst case as a boundary point of the interval under
investigation using the monotonicity of the criterion with respect to
$b$. This property,
which appears in many minimax- or maximin optimal design problems, has
been criticized by \citet{dette1997}.
This author recommends Bayesian or standardized maximin optimality
criteria, which reflect the different sizes of the optimality criteria
for different values of
$b$ in a more reasonable way. The determination of $T$-optimal
discriminating designs with respect to these criteria is substantially
harder and a challenging problem for future research.

\begin{appendix}\label{app}
\section*{\texorpdfstring{Appendix: Proof of identities (\lowercase{\protect\ref{weight1}}) and (\lowercase{\protect\ref{weight2}})}
{Appendix: Proof of identities (3.7) and (3.11)}}

Note that the identities in (\ref{weight1}) and (\ref{weight2}) can be
written in the form
%
\setcounter{equation}{0}
\begin{equation} \label{eq2}
\sum^n_{i=1}t^{*k}_i(b)(-1)^i \omega^*_i =0,\qquad k=0,1,\ldots,n-2,
\end{equation}
where $t^*_i (0)=\cos(\frac{i \pi}{n})=x_i$.
We will prove that these equalities hold for any real number $b$.
Since
%
\begin{equation} \label{eq3}\quad
t^{*k}_i(b)=\sum^k_{j=0}a_j\cos\biggl(\frac{ji\pi}{n} \biggr),\qquad i
=0,1,\ldots,n,k=0,1,\ldots,n-2,
\end{equation}
for some coefficients $a_j=a_j(b)$ $(j=0,1,\ldots, k)$ the identities in
(\ref{eq2}) follow from
%
\begin{equation} \label{eq4}
\sum^n_{i=1}(-1)^i\cos\biggl(\frac{ki\pi}{n} \biggr) \omega^*_i
=0,\qquad
k=0,1,\ldots,n-2.
\end{equation}
In order to prove (\ref{eq4}), consider first the case $k=0, n=2s$ for
some $s$, where the left-hand side of (\ref{eq4}) reduces to
\begin{eqnarray*}
\sum_{i=1}^n
\ww_i^*(-1)^i&=&\frac{1}{n}\Biggl[\sum_{i=1}^{s-1}\biggl[\biggl(1 -
\cos\biggl(\frac{i\pi}{n}\biggr)\biggr)(-1)^i+\biggl(1 +
\cos\biggl(\frac{i\pi}{n}\biggr)\biggr)(-1)^i\biggr] \\
&&\hspace*{195.6pt}{} +(-1)^s+1\Biggr]\\
&=&\frac{1}{n}\Biggl[\sum_{i=1}^{s-1}2(-1)^i +(-1)^s+1\Biggr]=0,
\end{eqnarray*}
which proves (\ref{eq4}).
If $k=0, n=2s+1$, we get
\begin{eqnarray*}
\sum_{i=1}^n \ww_i^*(-1)^i&=&\frac{1}{n}\Biggl[\sum_{i=1}^{s}\biggl[\biggl(1
- \cos\biggl(\frac{i\pi}{n}\biggr)\biggr)(-1)^i-\biggl(1 +
\cos\biggl(\frac{i\pi}{n}\biggr)\biggr)(-1)^i\biggr]+(-1)\Biggr]\\
&=&\frac{1}{n}\Biggl[2\sum_{i=1}^{s}\cos\biggl(\frac{i\pi
}{n}
\biggr)(-1)^{i+1}-1\Biggr] \\
&=& \frac{1}{n} \biggl[ 1- \frac{\cos[ {\pi
(1+2(n+1)s)}/({2n}) ]}{\cos( {\pi}/({2n}))} -1
\biggr]\\
&=& - \frac{1}{n} \frac{\cos( {(2s+1)\pi}/{2}
)}{\cos
({\pi}/({2n}))} = 0,
\end{eqnarray*}
where the third identity follows by standard results for
trigonometrical summation [see, e.g., \citet{jolley1961}, formula (428)].
This proves (\ref{eq4}) for the case $k=0, n=2s+1$.
Now consider the case of even $n, n=2s$ for some odd~$s$, $s=2l-1$ and
$k$ of the form $k=2(2r-1)$.
In this case
the left-hand side of (\ref{eq4}) reduces to
\begin{eqnarray*}
&& \frac{1}{n}\Biggl[\sum_{i=1}^{s-1}\biggl[\biggl(1 -
\cos\biggl(\frac{i\pi}{n}\biggr)\biggr)+\biggl(1 +
\cos\biggl(\frac{i\pi}{n}\biggr)\biggr)\biggr](-1)^i\cos\biggl(\frac
{ki\pi
}{n}\biggr) \\
&&\hspace*{143pt}{} +(-1)^s\cos\biggl(\frac{k\pi}{2}\biggr)+\cos
(k\pi
)\Biggr] \\
&&\qquad= \frac{1}{n}\Biggl[2\sum_{i=1}^{s-1}(-1)^i\cos\biggl(\frac{ki\pi
}{n}\biggr) +(-1)^s\cos\biggl(\frac{k\pi}{2}\biggr)+\cos(k\pi
)\Biggr] \\
& &\qquad= \frac{1}{n} \biggl\{ \biggl( \cos\biggl( \frac{k \pi}{4s}
\biggr)
\biggr)^{-1} \biggl[ \cos\biggl( \frac{\pi k}{4s} - \pi\biggr) +
\cos
\biggl( \frac{\pi k}{4s} + \frac{\pi}{2} (k+2s-2) \biggr) \biggr] +2
\biggr\} \\
&&\qquad= \frac{1}{n} \{ (-1) + (-1)^{2s-1}+2 \} = 0,
\end{eqnarray*}
where we have again used well-known results on trigonometric summation
[see \citet{jolley1961}, formula (428)].
Therefore we obtain equality (\ref{eq4}) in the case $n=2s, s=2l-1$ and
$k=2(2r-1)$. The other cases can be proved in a similar way, and the
details are omitted for the sake of brevity.
\end{appendix}

\section*{Acknowledgments}

The authors would like to thank Martina Stein, who typed numerous
versions of this manuscript with considerable technical expertise. The
authors would also like to thank two unknown referees for constructive
comments, which yield a substantial improvement of an earlier version
of this paper.


%

\printaddresses

\end{document}